\documentclass[reqno]{amsart}

\usepackage[english]{babel}
\usepackage{graphicx}
\usepackage{textcomp}
\usepackage{amsmath}
\newtheorem{thm}{ Theorem}

\newtheorem{proposition}{Proposition}%[section]
\newtheorem{definition}{Definition}%[section]

\begin{document}

\title [The Riemann problem for the pressureless gas dynamics] {The Riemann problem for the stochastically perturbed
non-viscous Burgers equation and the pressureless gas dynamics
model}

\author[Korshunova,\,Rozanova]{Anastasiya Korshunova, Olga Rozanova}

%\author[label1,label2]{}
\address{Mathematics and Mechanics Faculty, Moscow State University, Moscow
119992, Russia}

\thanks {Supported by Award DFG 436 RUS 113/823/0-1 and the special program of the
Ministry of Education of the Russian Federation "The development of
scientific potential of the Higher School", project 2.1.1/1399}

%\email[$^{1}$]{albeverio@uni-bonn.de}
\email{rozanova@mech.math.msu.su}



\date{\today}

\begin{abstract}
Proceeding from the method of stochastic perturbation of a Langevin
system associated with the non-viscous Burgers equation we construct
a solution to the Riemann problem for the non-interacting particles
and sticky particles systems. We analyze the difference in the
behavior of  discontinuous solutions for these two models and
relations between them.

%We find an interrelation  between a solution to of the
%stochastically perturbed non-viscous Burgers equation and
\end{abstract}

%\pacs{}
\maketitle

\medskip

%\end{abstract}

Let us consider the Cauchy problem for the non-viscous Burgers
equation:
\begin{equation}
\label{equ_Burg} \partial_t u+(u,\nabla)u=0,\qquad u(x,0)=u_0(x),
\end{equation}
%subject to initial data $u(x,0)=u_0(x)$,
where
$u(x,t)=(u_1,...,u_n)$ is a vector-function
$\mathbb{R}^{n+1}\rightarrow\mathbb{R}^n$.

It is well known that on the smooth solutions of the Burgers
equation this equation is equivalent to the system of ODEs
\begin{equation}
\label{ODU} \dot{x}(t)=u(t,x(t)),\quad \dot{u}(t,x(t))=0
\end{equation}
for the characteristics $x=x(t)$.

We associate with (\ref{ODU}) the following  system of stochastic
differential equations:
$$dX_k(t)=U_k(t)dt+\sigma d(W_k)_t,$$
\begin{equation}
\label{SDU} dU_k(t)=0,\quad k=1,...,n,
\end{equation}
$$X(0)=x,\quad U(0)=u,$$
where $(X(t),U(t))$ runs in the phase space
$\mathbb{R}^n\times\mathbb{R}^n,$  $\sigma>0$ is a constant and
$(W)_t=(W_k)_t$, $k=1,...,n,$ is the n - dimensional Brownian
motion.

Let us introduce a function
\begin{equation}
\label{u_sdu}
\hat{u}(t,x)=\dfrac{\int\limits_{\mathbb{R}^n}uP(t,x,u)du}{\int\limits_{\mathbb{R}^n}P(t,x,u)du},
\end{equation}
where $P(t,x,u)$ is the probability density in position and velocity
space. This value has a sense of the conditional expectation of $U$
for fixed position $X$  \cite{Chorin}. If we choose
\begin{equation}
\label{P0} P_0(x,u)=\delta(u-u_0(x))f_0(x)
=\prod\limits_{k=1}^n\delta(u_k-(u_0(x))_k)f_0(x),
\end{equation}
where $f_0(x)$ is an arbitrary sufficiently regular nonnegative
function such that $\,\int\limits_{{\mathbb R}^n}f_0(x)dx=1$, then
$\hat{u}(0,x)=u_0(x)$. Certain properties of $\hat{u}(t,x)$ was
studied in \cite{Roz1} and \cite{Roz2}.

The density $P=P(t,x,u)$ obeys the Fokker-Planck equation
\begin{equation}
\label{Fok-Plank} \dfrac{\partial P}{\partial
t}=\left[-\sum\limits_{k=1}^n u_k\dfrac{\partial}{\partial
x_k}+\sum\limits_{k=1}^n\dfrac12\sigma^2\dfrac{\partial^2}{\partial
x_k^2}\right] P,
\end{equation}
subject to initial data (\ref{P0}).

We apply the Fourier transform in the variables $x$ and $u$ to
(\ref{Fok-Plank}) and (\ref{P0}) and obtain the Cauchy problem for
$\tilde{P}=\tilde{P}(t,\lambda,\xi)$:
\begin{equation}
\label{preobr_Fok-Plank} \dfrac{\partial \tilde{P}}{\partial
t}=-\dfrac12\sigma^2|\lambda|^2\tilde{P}+\lambda\dfrac{\partial\tilde{P}}{\partial
\xi},
\end{equation}
\begin{equation}
\label{preobr_P0}\tilde{P}(0,\lambda,\xi)=\int\limits_{\mathbb{R}^n}e^{-i(\lambda,s)}e^{-i(\xi,u_0(s))}f_0(s)ds,
\end{equation}
which solution is given by the following formula:
\begin{equation}
\label{preobr_P}\tilde{P}(t,\lambda,\xi)=\tilde{P}(0,\lambda,\xi+\lambda
t)e^{-\frac12\sigma^2|\lambda|^2t}.
\end{equation}

The inverse Fourier transform allows to find the density
$P(t,x,u)$:
$$P(t,x,u)=\dfrac1{(2\pi)^{2n}}\int\limits_{\mathbb{R}^n}\int\limits_{\mathbb{R}^n}e^{i(\lambda,x)}e^{i(\xi,u)}\tilde{P}d\lambda d\xi=$$
\begin{equation}
\label{s_plotn}=\dfrac1{(\sqrt{2\pi
t}\sigma)^n}\int\limits_{\mathbb{R}^n}\,\delta(u-u_0(s))\,f_0(s)\,e^{-\frac{|u_0(s)t+s-x|^2}{2\sigma^2t}}ds.
\end{equation}
Then we substitute $P(t,x,u)$ in (\ref{u_sdu}) and get the
following expression for $\hat{u}(t,x)$:
\begin{equation}
\label{sol_u_sdu}
\hat{u}(t,x)=\dfrac{\int\limits_{\mathbb{R}^n}u_0(s)f_0(s)e^{-\frac{|u_0(s)t+s-x|^2}{2\sigma^2t}}ds}{\int\limits_{\mathbb{R}^n}f_0(s)e^{-\frac{|u_0(s)t+s-x|^2}{2\sigma^2t}}ds}.
\end{equation}
It should be noted that integrals in (\ref{sol_u_sdu}) are defined
also for a wider class of $f_0(x)$ then the probability density of
the particle positions in the space at the initial moment of time.
If the integral $\int\limits_{\mathbb{R}^n}f_0(x)dx$ diverges (for
example, for $f_0(x)=\rm const$), we consider the domain $[-L,L]^n$,
where $L>0$ and use another definition of $\hat{u}(t,x)$:
\begin{equation}\label{sol_u_sduL}\hat{u}(t,x)=\lim\limits_{L\rightarrow
+\infty}\dfrac{\int\limits_
{[-L,L]^n}u_0(s)f_0(s)e^{-\frac{|u_0(s)t+s-x|^2}{2\sigma^2t}}ds}
{\int\limits_{[-L,L]^n}f_0(s)e^{-\frac{|u_0(s)t+s-x|^2}{2\sigma^2t}}ds}.
\end{equation}
Evidently, this definition  coincides with (\ref{sol_u_sdu}) for
$f_0(x)\in L_1(\mathbb{R}^n)$.

The following property of $\hat{u}(t,x)$ holds:

\begin{proposition}
Let $u_0(x)$ and $f_0(x)>0$ be  bounded functions of class $C^1$
and the solution to the respective Cauchy problem (\ref{equ_Burg})
keeps smoothness for $t<t_*\le+\infty$. Then $\hat{u}(t,x)$ tends to
solution of problem (\ref{equ_Burg}) as $\sigma\rightarrow 0$ for
any fixed $(t,x)\in\mathbb{R}^{n+1},\,t<t_*$.
\end{proposition}

\proof Let us denote $J(u_0(x))$ the Jacobian matrix of the map
$\,x\longmapsto u_0(x).$ As it was shown in \cite{protter}\,(Theorem
1), if $J(u_0(x))$ has at least one  eigenvalue negative for a
certain point $x\in{\mathbb R}^n,$  then the classical solution to
(\ref{equ_Burg}) fails to exist beyond a positive time $t_*.$
Otherwise,  $t_*=\infty.$ The matrix $C(t,x) = (I+t J(u_0(x))),$
where $\,I\,$ is the identity matrix, fails to be invertible for
$t=t_*.$
%Further on we shall use the notation
%$(u_0)'_x$ instead of $|J(u_0(x)|.$

The formula (\ref{sol_u_sdu}) (or (\ref{sol_u_sduL})) implies
$$\lim\limits_{\sigma\rightarrow
0}\hat{u}(t,x)\,=\,\dfrac{\int\limits_{\mathbb{R}^n}u_0(s)f_0(s)\lim\limits_{\sigma\rightarrow
0}\frac1{(\sqrt{2\pi
t}\sigma)^n}e^{-\frac{|u_0(s)t+s-x|^2}{2\sigma^2t}}ds}{\int\limits_{\mathbb{R}^n}f_0(s)\lim\limits_{\sigma\rightarrow
0}\frac1{(\sqrt{2\pi
t}\sigma)^n}e^{-\frac{|u_0(s)t+s-x|^2}{2\sigma^2t}}ds}.$$

If $p(t,x,s) = u_0(s)t+s-x\,,$ we can use locally the implicit
function theorem and find $s=s(t,x,p).$ Moreover, $dp=\det
C(t,s)\,ds.$  Therefore,
$$\lim\limits_{\sigma\rightarrow
0}\hat{u}(t,x)\,=\,
\dfrac{\int\limits_{\mathbb{R}^n}u_0(s(p))f_0(s(p))\,\det
(C(t,s(p)))^{-1}\,\delta(p)\,dp}
{\int\limits_{\mathbb{R}^n}f_0(s(p))\,\det(C(t,s(p)))^{-1}\,\delta(p)\,dp\,}=\,u_0(s_0(t,x)),
$$
where we denote by $s_0(t,x)=s(t,x,0)$ the vector-function which
obeys the following vectorial equation:
\begin{equation}
\label{usl}u_0(s_0(t,x))t+s_0(t,x)-x=0.
\end{equation}
%Thus,
%$$\lim\limits_{\sigma\rightarrow
%0}\hat{u}(t,x)=u_0(s_0(t,x)).$$
Let us show that $u(t,x)=u_0(s_0(t,x))$ satisfies the Burgers
equation, that is
\begin{equation}\label{Burgsubs}\sum\limits_{j=1}^n\,\partial_j
(u_{0,i})(s_{0,j})_t\,+\,\sum\limits_{j,k=1}^n  u_{0,j}
\partial_k(u_{0,i}){s_{0,k}}_{x_j}=0, \quad i=1,...,n,
\end{equation}
and $u_0(s_0(0,x))=u_0(x)$. Here we denote by $u_{0,i}$ and
$s_{0,i}$ the $i$ - th components of vectors $u_0$ and $s_0,$
respectively.

We differentiate (\ref{usl}) with respect to $t$ and $x_j$ to get
the matrix equations: $$\sum\limits_{j=1}^n\,C_{ij}\,(s_{0,j})_t
+u_{0,i}=0,\, \quad i=1,...,n,$$ and
$$\sum\limits_{k=1}^n\,C_{ik}\,(s_{0,k})_{x_j}+\delta_{ij}=0,
\, \quad i,j=1,...,n,$$ where $\delta_{ij}$ is the Kronecker symbol.
The equations imply
\begin{equation}\label{subst}
(s_{0,j})_t\,=\,-\,\sum\limits_{i=1}^n\,(C^{-1})_{ij}\,u_{0,i},\qquad
(s_{0,k})_{x_j}\,=\,-\,(C^{-1})_{jk},\,\quad j=1,..,n.
\end{equation}
 It remains now only to
substitute (\ref{subst}) into (\ref{Burgsubs}).

Further,  (\ref{usl}) implies $u_0(s_0(0,x))=u_0(x)$.  $\,\square$

\medskip

It is important to note that $s_0(t,x)$ is unique for all $t$ for
which the solution to the Burgers equation $u(t,x)$ is smooth.
\medskip
\medskip

\medskip
\medskip

Let us denote $\rho(t,x)=\int\limits_{\mathbb{R}^n}P(t,x,u)du$. From
(\ref{s_plotn}) we have
\begin{equation}
\label{plotn} \rho(t,x)=\dfrac1{(\sqrt{2\pi
t}\sigma)^n}\int\limits_{\mathbb{R}^n}f_0(s)e^{-\frac{|u_0(s)t+s-x|^2}{2\sigma^2t}}ds.
\end{equation}

 It can be readily checked
that in the one-dimensional case the functions $\rho(t,x)$ and
$\hat{u}(t,x)$ solve the following system:
\begin{equation}
\label{sist_obw1}\partial_t \rho+\partial_x
(\rho\hat{u})=\dfrac12\sigma^2\,\partial_{xx}\rho,
\end{equation}
\begin{equation}
\label{sist_obw2}
\partial_t (\rho\hat{u})+\partial_x (\rho\hat{u}^2)=
-\int\limits_{\mathbb{R}}(u-\hat{u})^2\partial_x
P\,du+\dfrac12\sigma^2\,\partial_{xx}(\rho\hat{u}).
\end{equation}

The equation (\ref{sist_obw1}) follows from the Fokker-Planck
equation (\ref{Fok-Plank}) directly.  To prove  (\ref{sist_obw1}) we
note that definitions of $\hat{u}(t,x)$ and $\rho(t,x)$ imply
$$\partial_t
(\rho\hat{u})=\partial_t\int\limits_{\mathbb{R}}uP(t,x,u)du=\int\limits_{\mathbb{R}}uP'_t(t,x,u)du$$
$$=-\int\limits_{\mathbb{R}}u^2P'_x(t,x,u)du+\dfrac12\sigma^2\int\limits_{\mathbb{R}}uP''_{xx}(t,x,u)du$$
$$=-\int\limits_{\mathbb{R}}u^2P'_x(t,x,u)du+\dfrac12\sigma^2\partial^2_{xx}(\rho\hat{u}).$$
Further,
$$\partial_x
(\rho\hat{u}^2)=\hat{u}^2\partial_x\rho+2\rho\hat{u}\partial_x\hat{u}$$
$$=2\hat{u}\partial_x(\rho\hat{u})-\hat{u}^2\partial_x\rho=\int\limits_{\mathbb{R}}2\hat{u}uP'_x(t,x,u)du-
\int\limits_{\mathbb{R}}\hat{u}^2P'_x(t,x,u)du.$$

A certain analog of system (\ref{sist_obw1}),(\ref{sist_obw2}) was
obtained in \cite{Derm}.

Let us denote $f(t,x)=\lim\limits_{\sigma\rightarrow 0}\rho(t,x)$.
Taking into account Proposition 1 and the Fokker-Plank equations, as
a limit $\sigma\rightarrow 0$ for smooth $f(t,x)$ and $u(t,x)$ we
obtain the system of pressureless gas dynamic (e.g.\cite{Shand}) in
any space dimensions:
$$\partial_t f+div_x(fu)=0,$$\vspace*{-0.6cm}
\begin{equation}
\label{sist_pred}
\end{equation}\vspace*{-0.6cm}
$$\partial_t (fu)+\nabla(fu \otimes u)=0.$$

However the formula (\ref{sol_u_sdu}) has  sense also for
discontinuous initial data $(f_0(x),u_0(x))$. For the sake of
simplicity we dwell on the one-dimensional case and consider the
following initial data:
\begin{equation}
\label{f0_razr} f_0(x)=f_1+f_2\theta(x-x_0),
\end{equation}
\begin{equation}
\label{u0_razr} u_0(x)=u_1+u_2\theta(x-x_0),
\end{equation}
where $\theta(x-x_0)$ is the Heaviside function. Without loss of
generality we assume $x_0=0$.

\begin{definition} We call the couple of functions $(f(t,x),u(t,x))$ the
generalized solution to the problem (\ref{sist_pred})
(\ref{f0_razr}), (\ref{u0_razr}) in the sense of free particles, if
for almost all $(t,x)\in \mathbb{R}_+\times \mathbb{R}^n$
$$f(t,x)=\lim\limits_{\varepsilon\rightarrow
0}(\lim\limits_{\sigma\rightarrow 0}\rho^{\varepsilon}(t,x)),\qquad
u(t,x)=\lim\limits_{\varepsilon\rightarrow
0}(\lim\limits_{\sigma\rightarrow 0}\hat{u}^{\varepsilon}(t,x)),$$
where $(\rho^{\varepsilon}(t,x),\hat{u}^{\varepsilon}(t,x))$ satisfy
 the system (\ref{sist_obw1}), (\ref{sist_obw2})
with initial data $(f_0^{\varepsilon}(x),u_0^{\varepsilon}(x))$ from
the class $C^1(\mathbb{R}^n)$ such that
$\lim\limits_{\varepsilon\rightarrow 0}f_0^{\varepsilon}(x)=f_0(x)$,
$\lim\limits_{\varepsilon\rightarrow 0}u_0^{\varepsilon}(x)=u_0(x)$
for almost all fixed $x\in\mathbb{R}^n$.
\end{definition}

\begin{proposition}
The solution $(f(t,x),u(t,x))$ to the problem (\ref{sist_pred}),
(\ref{f0_razr}), (\ref{u0_razr}) in the sense of free particles does
not depend of the choice of the couple
$(f_0^{\varepsilon}(x),u_0^{\varepsilon}(x))\in
C^1(\mathbb{R}^n)\cap C_b(\mathbb{R}^n)$.
\end{proposition}

\proof  Let us choose two couples of smoothed functions
$\,(f_{01}^{\varepsilon}(x),u_{01}^{\varepsilon}(x))\,$ and \\
$(f_{02}^{\varepsilon}(x),u_{02}^{\varepsilon}(x))$ such that
$$\lim\limits_{\varepsilon\rightarrow
0}f^{\varepsilon}_{01}(x)=\lim\limits_{\varepsilon\rightarrow
0}f^{\varepsilon}_{02}(x)=f_0(x),\qquad
\lim\limits_{\varepsilon\rightarrow 0}u^{\varepsilon}_{01}(x)
=\lim\limits_{\varepsilon\rightarrow
0}u^{\varepsilon}_{02}(x)=u_0(x)$$ for any fixed $x\in\mathbb{R}^n$.
Then  the  couple
$$(f_0^{\varepsilon}(x),u_0^{\varepsilon}(x))=
(f_{01}^{\varepsilon}(x)-f_{02}^{\varepsilon}(x),\,
u_{01}^{\varepsilon}(x) -u_{02}^{\varepsilon}(x))\in C^1({\mathbb
R})\cap C_b(\mathbb{R}^n)$$ can be considered as initial data for
the problem (\ref{sist_pred})-(\ref{u0_razr}). To prove the
proposition we have to show that the respective solution  is
identically zero.

Indeed, from (\ref{plotn}) we have
$$f(t,x)=\lim\limits_{\varepsilon\rightarrow
0}\left(\int\limits_{\mathbb{R}^n}f_0^{\varepsilon}(s)\delta(s-s_0(t,x))ds\right)=
\lim\limits_{\varepsilon\rightarrow
0}(f^\varepsilon_{01}(s_0(t,x))-f^\varepsilon_{02}(s_0(t,x)))=0.$$
Here $s_0(t,x)$ is a solution of (\ref{usl}) as before. Analogously
proceeding from (\ref{sol_u_sdu}), we prove that $u(t,x)\equiv 0$.
$\square$

Our purpose is to find relations between a stable solution to the
Riemann problem of system that can be obtained as a limit $\sigma
\to 0$ from (\ref{sist_obw1}), (\ref{sist_obw2}) (as we will show
below it is not necessarily looks like (\ref{sist_pred})!) and the
couple  $(f(t,x),u(t,x)).$

\textbf{The Riemann problem
%for uniform initial distribution of
%density
in 1D case.} According to Definition 1 we must consider the smoothed
initial data instead of (\ref{f0_razr}) and (\ref{u0_razr}).  As
follows from Proposition 1  we can choose any couple of smoothed
initial data we want. It will be convenient to consider the
piecewise linear approximation of initial data of the form
\begin{equation}
\label{f} f^{\varepsilon}_0(x)=\begin{cases}
f_1,&\text{$x\leq -\varepsilon$,}\\
\dfrac{f_2}{2\varepsilon}x+f_1+\dfrac{f_2}2,&\text{$-\varepsilon<x<\varepsilon$,}\\
f_1+f_2,&\text{$x\geq \varepsilon$,}\\
\end{cases}
\end{equation}
\begin{equation}
\label{u} u^{\varepsilon}_0(x)=\begin{cases}
u_1,&\text{$x\leq -\varepsilon$,}\\
\dfrac{u_2}{2\varepsilon}x+u_1+\dfrac{u_2}2,&\text{$-\varepsilon<x<\varepsilon$,}\\
u_1+u_2,&\text{$x\geq \varepsilon$,}\\
\end{cases}
\end{equation}
where $f_1$, $f_2$, $u_1$ and $u_2$ are constants.

Note that these functions can be pointwisely approximated by
functions from the class $C^1(\mathbb{R}^n)$.

From (\ref{plotn}) we can find the density
$\rho^{\varepsilon}(t,x)$ corresponding to the smoothed initial
data $(f_0^{\varepsilon}(x),u_0^{\varepsilon}(x))$:
\begin{equation}
\label{P1}\rho^{\varepsilon}(t,x)=f_1\Phi\left(\dfrac{C_-^{\varepsilon}}{\sigma\sqrt{t}}\right)+(f_1+f_2)\Phi\left(-\dfrac{C_+^{\varepsilon}}{\sigma\sqrt{t}}\right)+I_1^{\varepsilon,\sigma},
\end{equation}
where
$\Phi(\alpha)=\dfrac1{\sqrt{2\pi}}\int\limits_{-\infty}^{\alpha}e^{-\frac{x^2}{2}}d
x$ is the Gauss function, $C_-^{\varepsilon}=u_1t-x-\varepsilon$,
$C_+^{\varepsilon}=(u_1+u_2)t-x+\varepsilon$, and
\begin{equation}
\label{I1}I_1^{\varepsilon,\sigma}=F^{\varepsilon}(t,x)\left[\Phi\left(\dfrac{C_+^{\varepsilon}}{\sigma}\right)-\Phi\left(\dfrac{C_-^{\varepsilon}}{\sigma}\right)\right]+O\left(\sigma
e^{-\frac1{\sigma^2}}\right),
\end{equation}
The expression for $I_1^{\varepsilon,\sigma}$ can be written out,
however,  we are interested only in behavior of
$I_1^{\varepsilon,\sigma}$ as $\sigma\rightarrow 0$. It can be
calculated that $\lim\limits_{\varepsilon\rightarrow
0}F^{\varepsilon}(t,x)=0$.

 To find
$\hat{u}(t,x)$ we compute the numerator in formula
(\ref{sol_u_sdu}):
$$\dfrac1{\sqrt{2\pi t}\sigma}\int\limits_{\mathbb{R}}u_0^{\varepsilon}(s)f_0^{\varepsilon}(s)
e^{-\frac{|u_0^{\varepsilon}(s)t+s-x|^2}{2\sigma^2t}}ds$$
$$=u_1\rho^{\varepsilon}(t,x)+u_2(f_1+f_2)\Phi\left(-\dfrac{C_+^{\varepsilon}}{\sigma\sqrt{t}}\right)+I_2^{\varepsilon,\sigma},$$
where
\begin{equation}
\label{I2}I_2^{\varepsilon,\sigma}=N^{\varepsilon}(t,x)\left[\Phi\left(\dfrac{C_+^{\varepsilon}}
{\sigma\sqrt{t}}\right)-\Phi\left(\dfrac{C_-^{\varepsilon}}{\sigma\sqrt{t}}\right)\right]+O\left(\sigma^2+\sigma
e^{-\frac1{\sigma^2}}\right),
\end{equation}
where
$N^{\varepsilon}(t,x)=\left(\dfrac{u_2}{u_2t+2\varepsilon}(x-(u_1+\dfrac{u_2}2)t)
+\dfrac{u_2}2\right)F^{\varepsilon}(t,x)$ and
$\lim\limits_{\varepsilon\rightarrow 0}N^{\varepsilon}(t,x)=0$.

Thus, we have the following result:
\begin{equation}
\label{sol_sdu1}
\hat{u}^{\varepsilon}(t,x)=u_1
+\dfrac{u_2(f_1+f_2)\Phi\left(-\frac{C_+^{\varepsilon}}{\sigma\sqrt{t}}\right)+I_2^{\varepsilon,\sigma}}{f_1\Phi\left(\frac{C_-^{\varepsilon}}{\sigma\sqrt{t}}\right)+(f_1+f_2)\Phi\left(-\frac{C_+^{\varepsilon}}{\sigma\sqrt{t}}\right)+I_1^{\varepsilon,\sigma}},
\end{equation}
where $I_1^{\varepsilon,\sigma}$ and $I_2^{\varepsilon,\sigma}$ are
given by in (\ref{I1}) and (\ref{I2}), respectively. Note that
$\dfrac{C_{\pm}^\varepsilon}{\sigma}\rightarrow\pm\infty$ as
$\sigma\rightarrow 0$.

Now  we can find the generalized solution to the Riemann problem
as $$f(t,x)=\lim\limits_{\varepsilon\rightarrow
0}(\lim\limits_{\sigma\rightarrow
0}\rho_{\varepsilon}(t,x)),\qquad
u(t,x)=\lim\limits_{\varepsilon\rightarrow
0}(\lim\limits_{\sigma\rightarrow 0}\hat{u}_{\varepsilon}(t,x)).$$

Let us introduce the points
$\hat{x}^{\varepsilon}_1=u_1t-\varepsilon$ and
$\hat{x}^{\varepsilon}_2=(u_1+u_2)t+\varepsilon$. Their velocities
are $u_1$ and $u_1+u_2$, respectively.

We consider two cases:

{\bf 1.} $u_2>0$ (velocity of the point $\hat{x}^{\varepsilon}_2$ is
higher than velocity of the point $\hat{x}^{\varepsilon}_1$). We can
find $f^{\varepsilon}(t,x)=\lim\limits_{\sigma\rightarrow
0}\rho^{\varepsilon}(t,x)$ from (\ref{P1}). Let us note that this
formula contains $F^{\varepsilon}(t,x)$. It is easy to see that
$$\lim\limits_{\varepsilon\rightarrow
0}\hat{x}^{\varepsilon}_1=x-u_1t,\qquad
\lim\limits_{\varepsilon\rightarrow
0}\hat{x}^{\varepsilon}_2=x-(u_1+u_2)t$$ and
$\lim\limits_{\varepsilon\rightarrow 0}F^{\varepsilon}(t,x)=0.$
Thus,
$$
f(t,x)=\lim\limits_{\varepsilon\rightarrow
0}f^{\varepsilon}(t,x)=\begin{cases}
f_1,&\text{$x<u_1t$,}\\
\dfrac{f_1}2,&\text{$x=u_1t$,}\\
0,&\text{$u_1t<x<(u_1+u_2)t$,}\\
\dfrac{f_1+f_2}2,&\text{$x=(u_1+u_2)t$,}\\
f_1+f_2,&\text{$x>(u_1+u_2)t$.}\\
\end{cases}
$$
Further, from (\ref{sol_sdu1}) we find the solution of the gas
dynamic system with smooth initial data
$u^{\varepsilon}(t,x)=\lim\limits_{\sigma\rightarrow
0}\hat{u}_{\varepsilon}(t,x)$ as follows:
$$
u^{\varepsilon}(t,x)=\begin{cases}
u_1,&\text{$x<\hat{x}^{\varepsilon}_1$,}\\
u_1+\dfrac{N(\varepsilon)}{F(\varepsilon)},&\text{$\hat{x}^{\varepsilon}_1\leq x\leq \hat{x}^{\varepsilon}_2$,}\\
u_1+u_2,&\text{$x>\hat{x}^{\varepsilon}_2$.}\\
\end{cases}
$$
It can be shown that $$\lim\limits_{\varepsilon\rightarrow
0}\dfrac{N(\varepsilon)}{F(\varepsilon)}=\lim\limits_{\varepsilon\rightarrow
0}\left(\dfrac{u_2}2
+\dfrac{u_2}{u_2t+2\varepsilon}(x-(u_1+\dfrac{u_2}2)t)\right)=\dfrac{x}{t}-u_1.$$
Thus, we get the following solution:
\begin{equation}
\label{rarefaction} u(t,x)=\lim\limits_{\varepsilon\rightarrow
0}u^{\varepsilon}(t,x)=\begin{cases}
u_1,&\text{$x< u_1t$,}\\
\dfrac{x}{t},&\text{$u_1t\leq x\leq(u_1+u_2)t$,}\\
u_1+u_2,&\text{$x>(u_1+u_2)t$.}\\
\end{cases}
\end{equation}
 We can see that the velocity includes the rarefaction wave. It
is well known stable solution to the Burgers equation with initial
data (\ref{u0_razr}).

It is interesting to note that if we compute the limit in
$\varepsilon$ first we get the solution
$u(t,x)=u_1+u_2\theta(x-(u_1+\dfrac{u_2}2)t)$, which is unstable
with respect to  small perturbations.

{\bf 2.} $u_2<0$ (the velocity of $\hat{x}^{\varepsilon}_2$ is
higher than the velocity of $\hat{x}^{\varepsilon}_1$). From
(\ref{P1}) and (\ref{sol_sdu1}) we find as before:
$$
f(t,x)=\begin{cases}
f_1,&\text{$x<(u_1+u_2)t$,}\\
\dfrac{3f_1+f_2}2,&\text{$x=(u_1+u_2)t$,}\\
2f_1+f_2,&\text{$(u_1+u_2)t<x<u_1t$,}\\
\dfrac{3f_1+2f_2}2,&\text{$x=u_1t$,}\\
f_1+f_2,&\text{$x>u_1t$,}\\
\end{cases}
$$
$$
u(t,x)=\begin{cases}
u_1,&\text{$x<(u_1+u_2)t$,}\\
u_1+\dfrac{f_1+f_2}{2f_1+f_2}u_2,&\text{$(u_1+u_2)t\leq x\leq u_1t$,}\\
u_1+u_2,&\text{$x> u_1t$.}\\
\end{cases}
$$
{\it Remark}\,
 We can  consider in
this framework the singular Riemann problem with initial density
$f_0(x)=f_1+f_2\theta(x-x_0)+f_3\delta(x).$

{\bf The Hugoniot conditions and the spurious pressure.} As we have
been proved, if $f$ and $u$ are smooth, they solve the pressureless
gas dynamics system. Now we ask the question which system satisfy
the solution of this system with jumps in the sense of free
particles.

 The
system of conservation laws (\ref{sist_pred}) implies two Hugoniot
conditions that should be held on the jumps of the solution
\cite{Rogd}. This signifies the solution satisfies the system in the
sense of integral identities. If we denote by $D$ the velocity of
the jump and $[h(y)]=h(y+0)-h(y-0)$ the value of the jump, then the
continuity equation and the momentum conservation give $[f]D=[fu]$
and $[fu]D=[fu^2],$ respectively.

In the case $u_2>0$ the velocity is continuous, therefore the
Hugoniot conditions hold trivially.

 We should check these conditions for
the jumps in the case $u_2<0.$ An easy computation shows that the
first one is satisfied, however, the second one does not hold. To
understand the reason let us estimate the integral term in
(\ref{sist_obw2}) on the generalized solutions in the case $u_2<0:$
$$
\int\limits_{\mathbb R}\,(u-\hat u(t,x))^2\,P_x(t,x,u)\,du=$$
$$=
\frac{1}{\sqrt{2\pi t}\sigma}\,\int\limits_{\mathbb
R}\,f_0(s)(u_0(s)-\hat
u(t,x))^2\,\left(e^{-\frac{(u_0(s)t+s-x)^2}{2\sigma^2t}}\right)_x\,ds=
$$
$$
= -\,\frac{1}{\sqrt{2\pi t}\sigma}\,\int\limits_{\mathbb
R}\,f_0(s)((u_0(s)-u_{FP}(t,x)) +(u_{FP}(t,x)-\hat
u(t,x)))^2\,\left(e^{-\frac{(u_0(s)t+s-x)^2}{2\sigma^2t}}\right)_s\,ds=
$$
$$
= -\,\frac{1}{\sqrt{2\pi t}\sigma}\,\int\limits_{\mathbb
R}\,f_0(s)(u_0(s)-u_{FP}(t,x))^2\,\left(e^{-\frac{(u_0(s)t+s-x)^2}{2\sigma^2t}}\right)_s\,ds+
$$
$$
 -2\,(u_{FP}(t,x)-\hat
u(t,x))\,\frac{1}{\sqrt{2\pi t}\sigma}\,\int\limits_{\mathbb
R}\,f_0(s)(u_0(s)-u_{FP}(t,x))\,\left(e^{-\frac{(u_0(s)t+s-x)^2}{2\sigma^2t}}\right)_s\,ds+
$$
$$
 -\,(u_{FP}(t,x)-\hat
u(t,x))^2\,\frac{1}{\sqrt{2\pi t}\sigma}\,\int\limits_{\mathbb
R}\,f_0(s)\,\left(e^{-\frac{(u_0(s)t+s-x)^2}{2\sigma^2t}}\right)_s\,ds=
$$
$$I_1+I_2+I_3.$$
The integrals $I_2$ and $I_3$ tend to zero as $\sigma\to 0$ due to
properties of the Riemann data since $\hat u(t,x)\to u_{FP}(t,x)$
for almost all $x\in \mathbb R.$ Let us estimate $I_1.$
$$I_1\,
=-\,\frac{1}{\sqrt{2\pi t}\sigma}\,\int\limits_{u_2
t}^0\,f_1\,(u_1-\hat
u(t,x))^2\,\left(e^{-\frac{(u_1\,t+s-x)^2}{2\sigma^2t}}\right)_s\,ds\,-$$$$
-\,\frac{1}{\sqrt{2\pi t}\sigma}\,\int\limits_{0}^{-u_2
t}\,(f_1+f_2)\,((u_1+u_2)-\hat
u(t,x))^2\,\left(e^{-\frac{((u_1+u_2)\,t+s-x)^2}{2\sigma^2t}}\right)_s\,ds
=
$$
$$=
-\,\frac{1}{\sqrt{2\pi t}\sigma}\,\frac{f_1(f_1+f_2)^2
u_2^2}{(2f_1+f_2)^2}\left(e^{-\frac{(u_1\,t-x)^2}{2\sigma^2t}}-e^{-\frac{((u_1+u_2)\,t-x)^2}{2\sigma^2t}}\right)-$$
$$-\,\frac{1}{\sqrt{2\pi t}\sigma}\,
\frac{f_1^2(f_1+f_2)
u_2^2}{(2f_1+f_2)^2}\left(e^{-\frac{(u_1\,t-x)^2}{2\sigma^2t}}-e^{-\frac{((u_1+u_2)\,t-x)^2}{2\sigma^2t}}\right).$$
Thus,
$$I_1\,\to\,-\frac{f_1(f_1+f_2)
u_2^2}{(2f_1+f_2)}(\delta(x-(u_1+u_2)\,t)\,- \, \delta(x-u_1\,t)
),\quad \sigma\to 0.$$

 Thus, the integral term corresponds to a spurious pressure
between the  jumps $x=(u_1+u_2)\,t$  and $x=u_1\,t,$ namely,
\begin{equation}\label{press}p=\frac{f_1(f_1+f_2)
u_2^2}{(2f_1+f_2)}(\theta(x-(u_1+u_2)\,t)\,- \, \theta(x-u_1\,t)).
\end{equation}
The Hugoniot condition $[fu]D=[fu^2+p]$ is satisfied with this kind
of pressure.

Thus, we get the following theorem.

\begin{thm} The generalized solution to the Riemann problem  (\ref{f0_razr}), (\ref{u0_razr})  for the pressureless gas
dynamics system in the sense of free particles (according to
Definition 1) in the case of a discontinuous velocity ($u_2<0$)
solves in fact the gas dynamics system with a pressure defined by
(\ref{press}).
\end{thm}

 {\bf
Sticky particles model vs non-interacting particles model.} In our
model the particles are allowed to go through the discontinuity as
one particle does not feel the others. However, if we are in the
frame of the sticky particles model we should assume that the
particles meeting one other  stick together on the jump.  The
non-interacting particles model and the sticky particles model are
equivalent for smooth velocities, however, if the velocity has a
jump, the solutions behavior differs drastically. Nevertheless, we
can study the solution to the Riemann problem in the case of $u_2<0$
for the sticky particles model, too, basing on the solution obtained
in the present paper. Indeed, the jump position $x(t)$ is  a point
between $x_1(t)=(u_1+u_2)\,t$ and
$x_2(t)=u_1\,t.$ %Let us denote $D(t)=\dot x(t).$
%its
%position is $x-(u_1+\dfrac{u_2}2)t$ (the Huginiot condition for the
%Burgers equation written in the form of a conservation law
%$u_t+\left(\frac{u^2}{2}\right)$ gives the same result).
The mass $m(t)$ accumulates in the jump due to the impenetrability
of the discontinuity with the velocity
$$m(t)=(x(t)-(u_1+u_2)t)((2f_1+f_2)-f_1)+(u_1t-x(t))((2f_1+f_2)-(f_1+f_2))\,$$$$=
\left(-((u_1+u_2)(f_1+f_2)-u_1\,f_1)t+x(t) f_2\right)\,=
-[uf]t+[f]x(t),$$ where  $[\,\,]$ stands for a jump value.
%Therefore
%due to the initial conditions $x(0)=0,\,$ we have
%$$
%m(t)=-[uf]t+[f]x(t)+m(0)
%$$
Further, if we change heuristically the overlapped  mass between
$x_1$ and $x_2$ to the mass concentrated at a point, then from the
condition of equality of momenta in the both cases we can find the
velocity of the point singularity:
$$(u_1+u_2)((2f_1+f_2)-f_1)(x(t)-(u_1+u_2)t)+u_1\,f_1\,(u_1t-x(t))\,$$
$$=-[u^2 f]t+[u f] x(t)=m(t)\dot x(t).$$
%$$=D\left(-((u_1+u_2)(f_1+f_2)-u_1\,f_1)+Df_2\right).$$
Thus, to find the position of the point singularity we get the
equation
\begin{equation}
\label{jump} ([f]x(t)-[uf]t)\,\dot x(t)=[u f] x(t)-[u^2 f]t,
\end{equation}
subject to initial data $x(0)=0.$ The respective solution is
\begin{equation}\label{R_jumpposition}
x(t)=\frac{1}{[f]}\,\left([uf]\pm\sqrt{[uf]^2-[f][u^2f]}\right)\,t,
\end{equation}
where the sign should be chosen from the condition
$x_1(t)<x(t)<x_2(t).$ It can be readily shown that the latter
condition is satisfied either  for plus or minus in the formula
(\ref{R_jumpposition}). The condition coincides with the Lax
stability condition $u_1<\dot x(t)<u_1+u_2.$

The formulas describing  the amplitude of the delta-function in the
density component and the singularity position  obtained earlier in
\cite{yang},\cite{shelkv},\cite{danilov} give the same result.

It is worth mentioning that the spurious pressure (\ref{press}) does
not arise in the sticky particles model.

{\bf Back to the Burgers equation.} Now we are able to get the
solution to the Cauchy problem (\ref{equ_Burg}), (\ref{u0_razr}).
Let us assume that $f_0(x)=const,$ that is $f_2=0.$ Then for the
case $u_2<0$ from (\ref{jump}) we obtain the well known formula for
a position of the jump:
$$x(t)=\frac{[u^2]}{2[u]}\,t=\frac{u_++u_-}{2}.$$
If $u_2>0,$ the solution is continuous and it is given by the
formula (\ref{rarefaction}).

\end{document}